\documentclass[12pt]{article}

\usepackage{amsfonts}
\usepackage{amsmath}
\usepackage{amssymb}
\usepackage{amscd}
\usepackage{amsthm}
\usepackage{indentfirst}

\newtheorem{thm}{Theorem}[section]

\newtheorem{prop}[thm]{Proposition}
\newtheorem{cor}[thm]{Corollary}

\newtheorem{exam}[thm]{Example}

\newtheorem{question}[thm]{Question}

\newcommand{\R}{{\mathbb{R}}}

\newcommand{\Z}{{\mathbb{Z}}}
\newcommand{\N}{{\mathbb{N}}}
\newcommand{\C}{{\mathbb{C}}}
\newcommand{\SP}{{\mathbb{S}}}

\newcommand{\cG}{{\mathcal{G}}}

\def\id{{1\hskip-2.5pt{\rm l}}}

\newcommand{\Ham}{{\hbox{\it Ham\,}}}

\newcommand{\Symp}{{\hbox{\it Symp} }}

\newcommand{\osc}{{\text{osc}}}

\newcommand{\Qed}{\hfill \qedsymbol \medskip}

\DeclareMathOperator{\sgrad}{sgrad}

\begin{document}

\title{Quasi-morphisms and the Poisson bracket \\
}

\renewcommand{\thefootnote}{\alph{footnote}}

\author{\textsc Michael Entov$^{a}$,\ Leonid
Polterovich$^{b}$\ and Frol Zapolsky$^{b}$}

\footnotetext[1]{Partially supported by E. and J. Bishop Research
Fund. } \footnotetext[2]{Partially supported by the Israel Science
Foundation grant $\#$ 11/03.}

\date{\today}

\maketitle

\centerline{{\it To Gregory Margulis with admiration}}

\begin{abstract}
\noindent For a class of symplectic manifolds, we introduce a
functional which assigns a real number to any pair of continuous
functions on the manifold. This functional has a number of
interesting properties. On the one hand, it is Lipschitz with
respect to the uniform norm. On the other hand, it serves as a
measure of non-commutativity of functions in the sense of the
Poisson bracket, the operation which involves first derivatives of
the functions. Furthermore, the same functional gives rise to a
non-trivial lower bound for the error of the simultaneous
measurement of a pair of non-commuting Hamiltonians. These results
manifest a link between the algebraic structure of the group of
Hamiltonian diffeomorphisms and the function theory on a
symplectic manifold. The above-mentioned functional comes from a
special homogeneous quasi-morphism on the universal cover of the
group, which is rooted in the Floer theory.
\end{abstract}

\vfil \eject


\renewcommand{\thefootnote}{\arabic{footnote}}

\section{Introduction and main results}
\label{sec-intro}

\subsection{Symplectic quasi-states}\label{subsec-qst}\noindent
Let $(M^{2n}, \omega)$ be a closed connected symplectic manifold.
A {\it symplectic quasi-state} is a functional $\zeta:
C^0(M)\to\R$ with the following properties:
\begin{itemize}
\item[{(i)}] $\zeta(1)=1$;
\item[{(ii)}] $F \leq G \Rightarrow \zeta(F)\leq \zeta(G)$;
\item[{(iii)}] $\zeta(aF+bG)=a\zeta(F)+b\zeta(G)$ for all $a,b \in \R$
and all functions $F,G \in C^{\infty}(M)$ whose Poisson bracket
$\{F,G\}$ vanishes.
\end{itemize}
In particular, $\zeta$ is a topological quasi-state in the sense
of Aarnes \cite{Aar}. It was shown in \cite{EP-qst} that certain
symplectic manifolds carry a {\bf non-linear} symplectic
quasi-state $\zeta$. In this case the functional
$$\Pi: C^0(M)\times C^0(M) \to \R$$ defined by
\begin{equation}\label{eq-pi-zeta}
\Pi(F,G) = |\zeta(F+G)-\zeta(F)-\zeta(G)|\;
\end{equation}
can be interpreted as a measure of Poisson non-commutativity of
functions $F$ and $G$. This functional lies in the center of our
study. In particular, we analyze the relation between $\Pi(F,G)$
and the Poisson bracket $\{F,G\}$. We show (see Theorem
\ref{thm-main-A} below) that for certain symplectic quasi-states
\begin{equation}\label{eq-intro1}
\Pi(F,G) \leq \text{const}\cdot \sqrt{||\{F,G\}||}
\end{equation}
for all $F,G \in C^{\infty}(M)$. Here  and below $||H||$ stands
for the uniform norm $\max_M |H|$. This inequality has a number of
applications.

One application deals with the following problem (cf.
\cite{Car-Vit}). The definition of the Poisson bracket $\{ F, G\}$
of two functions $F, G$ on a symplectic manifold $M$ involves
first derivatives of the functions. Thus \emph{a priori} there is
no restriction on possible changes of $\{ F, G\}$ when $F$ and $G$
are perturbed in the uniform norm. Note that axiom (ii) of the
quasi-state yields that $\Pi$ is Lipschitz in the uniform norm.
Therefore inequality \eqref{eq-intro1} gives rise to such a
restriction (see Corollary \ref{cor-robust} below).

As another application, we present  a restriction on partitions of
unity $\rho_1,\ldots,\rho_N$ on symplectic manifolds subordinate
to coverings by sufficiently small sets (see Theorems
\ref{thm-part} and \ref{thm-part-1} below). It turns out that
$$\max_{i,j} ||\{\rho_i,\rho_j\}|| \geq
\frac{\text{const}}{N^3}\;.$$

Interestingly enough, the functional $\Pi$ appears in the context
of simultaneous measurements of non-commuting observables $F,G$ in
classical mechanics (see Section \ref{sect-measurement-general}
below). We show that $\Pi(F,G)$ gives a lower bound for the error
of such a measurement.

The above-mentioned quasi-states are closely related to certain
homogeneous quasi-morphisms (that is, homomorphisms up to a
bounded error) on the universal cover $\cG$ of the group of
Hamiltonian diffeomorphisms of $M$. These quasi-morphisms, which
were found in \cite{EP-qmm}, are rooted in Floer homology. The
connection between quasi-states and quasi-morphisms is crucial for
our methods.

\subsection{Preliminaries on Hamiltonian diffeomorphisms}
\noindent

In what follows we normalize the symplectic form $\omega$ on
$M^{2n}$ so that the symplectic volume $ \int_M \omega^n$ equals
1. Recall that symplectic manifolds appear as phase spaces in
classical mechanics. An important principle of classical mechanics
is that {\it the energy of a system determines its evolution}. The
energy (or {\it Hamiltonian function}) $F_t(x):=F(x,t)$ is a
smooth function on $M \times [0;1]$. Here $t$ is the time
coordinate. Define the time-dependent Hamiltonian vector field
$\sgrad F_t$ by the point-wise linear equation $i_{\sgrad {F_t}}
\omega = -dF_t$. The evolution of the system is described by the
flow $f_t$ on $M$ generated by the Hamiltonian vector field
$\sgrad F_t$. We shall call the time-one-map $f_1$ of this flow
{\it a Hamiltonian diffeomorphism}. Hamiltonian diffeomorphisms
form a group which is denoted by $\Ham (M,\omega)$. The universal
cover $\cG$ of this group plays an important role in this paper.
Elements of $\cG$ are smooth paths in $\Ham (M,\omega)$ based at
the identity, considered up to homotopy with fixed end-points. We
denote by $\phi_F$ the element of $\cG$ represented by the path
$\{f_t\}_{t \in[0;1]}$ and refer to it as the element of $\cG$
generated by $F$.

It is instructive to view $\cG$ as a (infinite-dimensional) Lie
group whose Lie algebra is naturally identified with the space
$C_0^{\infty}(M)$ of smooth functions with zero mean on $M$. The
role of the Lie bracket is played by the Poisson bracket which is
defined by
$$\{F,G\} = \omega(\sgrad G,\sgrad F) = dF (\sgrad G) = - dG (\sgrad F) \;.$$
In the canonical local coordinates $(p,q)$ where $\omega =
dp\wedge dq$ the Poisson bracket is written as
$$\{F,G\} = \frac{\partial F}{\partial q}\cdot\frac{\partial G}{\partial p}-
\frac{\partial F}{\partial p}\cdot\frac{\partial G}{\partial
q}\;.$$

Recall that for a function $F \in C^0(M)$ we denote by $||F||$ its uniform norm
$\max_M | F | $ and by $\langle F \rangle$ its mean value $\int_M
F\omega^n$. A Hamiltonian function $F\in C^{\infty} (M \times
[0;1])$ is called {\it normalized} if $\langle F_t \rangle = 0$
for all $t$.

We refer to \cite{MS,Pbook} for further preliminaries on
Hamiltonian diffeomorphisms.

\subsection{Quasi-morphisms}
\noindent

A real-valued function $\mu$ on a group $\Gamma$ is called a {\it
homogeneous quasi-morphism} \cite{Bav} if
\begin{itemize}
\item[{(i)}] There
exists $C>0$ so that
\[
|\mu (\varphi\psi) - \mu (\varphi) - \mu(\psi)|\leq C \ \ {\rm
for\;\; all\;\; elements}\ \varphi,\psi\in\Gamma;
\]
\item[{(ii)}] $\mu (\varphi^m) = m\mu (\varphi)$ for each $\varphi \in \Gamma$ and each
$m\in\Z$.
\end{itemize}
The minimal constant $C$ in the above inequality is called {\it
the defect} of $\mu$.

\medskip

In this paper we will deal with homogeneous quasi-morphisms on the
group $\cG$ with the following property:
\begin{equation}\label{eq-Lipschitz}
\int_0^1 \min_M (F_t-G_t)\;dt \leq \mu(\phi_G)-\mu(\phi_F)
\leq \int_0^1 \max_M(F_t-G_t)\;dt
\end{equation}
for all {\bf normalized} Hamiltonians $F,G \in C^{\infty} (M
\times [0;1])$. We call them {\it stable} quasi-morphisms.

\medskip
\noindent \begin{exam}\label{ex-mfds}{\rm The group $\cG$ is known
to carry a stable homogeneous quasi-morphism for the following
list of symplectic manifolds \cite{EP-qmm,EP-semisimple,Ostr}:
complex projective spaces and Grassmannians; $\C P^2$ blown up at
$k\leq 3$ points with a symplectic form in a rational cohomology
class; strongly semi-positive (see Section~\ref{subsec-part})
direct products of the above-mentioned manifolds. The existence of
stable homogeneous quasi-morphisms is related to the algebraic
structure of the quantum homology of $(M,\omega)$. See Section
\ref{subsec-stab} below for more discussion on the stability
property.}
\end{exam}

\medskip
\noindent

Given a stable homogeneous quasi-morphism $\mu$, define a
functional $\zeta: C^{\infty}(M) \to \R$ by
\begin{equation}\label{eq-qmm-qst}
\zeta(F) = \int_M F\omega^n -\mu(\phi_F)\;.
\end{equation}
Recall that the Lie algebra of the group $\cG$ can be identified
with the space $C^{\infty}_0(M)$ of smooth functions on $M$ with
zero mean. With this language the restriction of $\zeta$ to
$C^{\infty}_0(M)$ is simply the pullback of quasi-morphism $-\mu$
on the group to the Lie algebra via the exponential map. One can
show that $\zeta$ satisfies the axioms of a symplectic quasi-state
listed in Section \ref{subsec-qst}: Axiom (i) is obvious (since,
according to our normalization, $\int_M \omega^n = 1$), axiom (ii)
is a simple corollary of the stability property
\eqref{eq-Lipschitz} of $\mu$ (see Section \ref{subseq-pr-eqlip}
below) and axiom (iii) follows from the fact (which is an easy
exercise) that the restriction of any homogeneous quasi-morphism
to an abelian subgroup is a homomorphism (see \cite{EP-qst}).

As an immediate consequence of axiom (ii) we get that $\zeta$ is
1-Lipschitz with respect to the uniform norm and thus extends to
the space of continuous functions $C^0(M)$. Furthermore, the
functional $$\Pi(F,G) = |\zeta(F+G)-\zeta(F)-\zeta(G)|$$ (see
formula \eqref{eq-pi-zeta} above) is Lipschitz as well:
\begin{equation}\label{eq-Lip-Pi}
|\Pi(F,G)-\Pi(F',G')|\leq 2(||F-F'||+||G-G'||)
\end{equation}
for all functions $F,G,F',G' \in C^{\infty}(M)$.

It is important to emphasize that in the setting of Example
\ref{ex-mfds} above the quasi-state $\zeta$ is non-linear, that is
$\Pi(F,G) > 0$ for some $F,G \in C^{\infty}(M)$.

\medskip
\noindent
\begin{exam}\label{ex-qst-sphere}{\rm Even for the 2-sphere with
the standard area form of total area $1$ the explicit calculation
of $\mu(f)$ for Hamiltonian diffeomorphisms $f$ generated by a
generic time-dependent Hamiltonian is a transcendentally difficult
problem. However the corresponding quasi-state $\zeta$ has an easy-to-handle
combinatorial interpretation, see \cite{EP-qmm}: Since
$\zeta$ is Lipschitz in the uniform norm, it suffices to define
its value on the dense subset consisting of Morse functions $F$ on
the sphere with distinct critical values. Look at the set of
connected components of the level lines of $F$. One can show that
there exists unique component, say, $\gamma_F$ with the following
property: the area of every connected component of the complement
$\SP^2 \setminus \gamma_F$ is $\leq \frac{1}{2}$. Then $\zeta(F)$
is simply the value $F(\gamma_F)$. Moreover, if $u:\R\to \R$ is
any smooth function, $\zeta(u\circ F) = u(\zeta(F))$.}
\end{exam}

\begin{exam}\label{ex-calc-sphere}
{\rm  Think of the 2-sphere $\SP^2$ as of the Euclidean unit
sphere in $\R^3(x,y,z)$ with the center at zero. Let $\omega$ be
the induced area form on $\SP^2$ divided by $4\pi$. We claim that
$\Pi(x^2,y^2) =1$ and hence $\zeta$ is non-linear. Indeed, apply
the explicit formula for $\zeta$ presented in the example above:
Note that for $F(x,y,z)=x$ the component $\gamma_F$ is the equator
$\{x=0\}$, and so $\zeta(x^2)=0$. Similarly
$\zeta(y^2)=\zeta(z^2)=0$. On the other hand
$$\zeta(x^2+y^2)= \zeta(1-z^2)\stackrel * = 1-\zeta(z^2)= 1\;,$$
where equality $(*)$ is valid in view of axioms (i),(iii) of a
quasi-state. Summing up,
$$\Pi(x^2,y^2) = |\zeta(x^2+y^2)-\zeta(x^2)-\zeta(y^2)| = 1\;,$$
which proves the claim. }
\end{exam}

\subsection{A lower bound for the Poisson bracket}\label{subsec-Pois}\noindent
Let $(M,\omega)$ be a closed symplectic manifold. Assume that the
group $\cG$ admits a stable homogeneous quasi-morphism with
defect $C$. Let $\zeta$ be the corresponding quasi-state and let
$\Pi$ be the functional defined by \eqref{eq-pi-zeta}.

\begin{thm}\label{thm-main-A}
$$ \Pi(F,G) \leq \sqrt{2C \cdot ||\{F,G\}||}$$
for all $F,G \in C^{\infty}(M)$.
\end{thm}

\medskip
\noindent The proof is given in Section \ref{subsec-proof-main}.

Let us describe an application of this result to $C^0$-robust
lower bounds on the Poisson bracket. We start our discussion with
the following result from \cite{Car-Vit}:
\begin{equation}\label{eq-vit}F_n \stackrel{C^0}{\longrightarrow} F,\; G_n
\stackrel{C^0}{\longrightarrow} G,\; \{ F_n, G_n\}
\stackrel{C^0}{\longrightarrow} 0\;\Rightarrow \{ F, G\} \equiv
0\end{equation} (all functions are assumed to be smooth). For the
sake of completeness, we present a proof in Section
\ref{subsec-vit} below. In particular, if $\{F,G\} \not\equiv 0$,
$$\exists\, \epsilon_0 = \epsilon_0(F,G)>0 \; : \; \forall\,
\epsilon < \epsilon_0 \;\; \exists\, \delta=\delta(\epsilon,F,G)
> 0$$ so that for all smooth functions $F',G'$ with
$$||F - F'|| + ||G - G'|| \leq \epsilon$$
we have $||\{F^\prime, G^\prime\}|| \geq \delta$.

This gives rise to the following definitions. Given $\epsilon > 0$
consider an open $C^0$-neighborhood $U_\epsilon$ of $(F, G)$ in
$C^\infty (M) \times C^\infty (M)$ defined as
\[
U_\epsilon := \{ \ (F^\prime, G^\prime ) \in  C^\infty (M) \times C^\infty (M)
\  \ | \ \   || F^\prime - F|| + || G^\prime - G|| < \epsilon\ \}.
\]
Set
\[
\Upsilon_{F, G} (\epsilon) := \inf_{U_\epsilon}
\| \{ F^\prime, G^\prime  \} \|
\]
and
\[
\Upsilon (F, G) := \lim_{\epsilon\searrow 0} \Upsilon_{F,G}
(\epsilon)= \liminf_{F',G'\stackrel{C^0}{\longrightarrow}F,G}
||\{F',G'\}||\;.
\]

 The largest $\epsilon_0$ as above, denoted by
$\epsilon_{max} (F, G)$, can be represented as follows:
\[
\epsilon_{max} (F, G) = \sup \ \{\  \epsilon \ | \   \Upsilon_{F,G} (\epsilon) >
0\
\}.
\]
It reflects the size of the "maximal" neighborhood $U_\epsilon$ of
$(F, G)$ which does not contain a pair of Poisson-commuting
functions. Given a positive $\epsilon < \epsilon_{max}$ one can
pick the corresponding $\delta(\epsilon,F,G)$ (see above) as
$\Upsilon_{F,G} (\epsilon)$.

It is an interesting problem to find explicit (lower) estimates for
the numbers $\Upsilon (F, G)$, $\epsilon_{max} (F,G)$ and the
function $\Upsilon_{F,G}: (0, \epsilon_0) \to \R$ (for at least some
$\epsilon_0 \in (0, \epsilon_{max})$) in terms of $F$ and $G$. As we
shall see in Section \ref{subsec-vit} below, the proof of
\eqref{eq-vit} leads to such estimates which involve Hofer's norm of
the commutator of the Hamiltonian diffeomorphisms generated by $F$
and $G$ (see formulae \eqref{eq-hof-final}, \eqref{eq-hof-4} below).
Theorem \ref{thm-main-A} gives us an expression of a different
nature, namely in terms of $\Pi(F,G)$, provided $\Pi(F,G)\neq 0$. As
a consequence, in some examples, explicit estimates on $\Upsilon (F,
G)$, $\epsilon_{max} (F,G)$ and $\Upsilon_{F,G} (\epsilon)$ can be
easily obtained using the machinery of symplectic quasi-states.

\begin{cor}
\label{cor-robust}
Let $F,G \in C^{\infty}(M)$ be two functions with $\Pi(F,G) \neq
0$. Then
\[
\| \{ F^\prime, G^\prime \}\| \geq \frac{ (\Pi (F, G) - 2\| F -
F^\prime\| - 2 \| G - G^\prime \| )^2}{2C}
\]
for all $F',G' \in C^{\infty}(M)$ with
$$||F-F'||+||G-G'|| \leq \frac{\Pi(F,G)}{2}\;.$$
In particular, $\Upsilon (F, G)\geq (\Pi (F, G))^2/2C$,
$\epsilon_{max} (F, G) \geq \Pi(F,G)/2$ and  \\ $\Upsilon_{F,G}
(\epsilon) \geq (\Pi(F,G)-2\epsilon)^2/2C$ for all $\epsilon\in
(0, \Pi(F,G)/2)$.

\end{cor}

\medskip
\noindent This is an immediate consequence of Theorem~\ref{thm-main-A} and
inequality \eqref{eq-Lip-Pi}.

The calculation of the defect $C$ is so far an open problem even in
the simplest examples. However, upper bounds $C \leq C_0$ are
available. In view of Corollary \ref{cor-robust}, if $\Pi(F,G)\neq
0$ we can pick $\epsilon_0 (F,G)= \Pi(F,G)/2$ and give the following estimates:
\begin{equation}\label{eq-epdelta}
\Upsilon (F, G) \geq (\Pi (F, G))^2/2C_0,\ \ \Upsilon_{F,G} (\epsilon) \geq
\frac{(\Pi(F,G)-2\epsilon)^2}{2C_0}\ \ \forall\, \epsilon\in (0, \epsilon_0) \;.
\end{equation}

Let us illustrate these inequalities in specific examples. We start
with the functions $F=x^2$ and $G=y^2$ on the two-sphere $\SP^2$
(see Examples \ref{ex-qst-sphere} and \ref{ex-calc-sphere} above).

\begin{prop}\label{prop-illustr} The quasi-state $\zeta$ given in
Example \ref{ex-qst-sphere} is induced by a stable homogeneous
quasi-morphism with defect $C \leq 2$.
\end{prop}

\medskip
\noindent This is proved in \cite{EP-qmm} with the exception of
the upper bound on the defect. We derive this bound in Section
\ref{sec-defect} below.

Taking into account that $\Pi(x^2,y^2)=1$ (see Example
\ref{ex-calc-sphere} above) and that $C_0= 2$ we get that

\begin{equation}\label{eq-illustr-1}
\epsilon_{max} (x^2,y^2) \geq 1/2 \;\;\text{and}\;\;\Upsilon_{x^2,
y^2} (\epsilon) \geq (1-2\epsilon)^2/4\;\; \forall \epsilon\in (0,
1/2)\;.
\end{equation}

\medskip
\noindent We refer to Section \ref{sec-prob} for further discussion
of inequality \eqref{eq-illustr-1}.

\medskip Let us now outline what happens in a higher dimensional example.
Consider the product of two spheres $M:= \SP^2 \times \SP^2$
equipped with the split symplectic form $\omega \oplus \omega$. Let
$\mu$ be the stable homogeneous  {\it Calabi} quasi-morphism on
$\cG$ defined in \cite{EP-qmm} (see Example \ref{ex-mfds} above).
The proof of the fact that $\mu$ is a quasi-morphism presented in
\cite{EP-qmm} is constructive, and hence gives rise to an upper
bound for the defect $C$ of $\mu$. In particular, unveiling the
argument presented \cite[Section 3.3]{EP-qmm} in the case of $\SP^2
\times \SP^2$ one gets the upper bound $C \leq 6$ (the details are
somewhat technical and will be omitted). Denote by $\zeta$ the
quasi-state associated to $\mu$ by formula \eqref{eq-qmm-qst}. Put
$F=x_1^2$ and $G=y_1^2$ where $(x_1,y_1,z_1)$ are the Euclidean
coordinate functions on the first factor of $M$. It is an immediate
consequence of \cite{BEP,EP-qst} that $\Pi(F,G)=1$. Thus Corollary
\ref{cor-robust} yields
$$\epsilon_{max} (x_1^2,y_1^2) \geq 1/2
\;\;\text{and}\;\;\Upsilon_{x_1^2, y_1^2} (\epsilon) \geq
(1-2\epsilon)^2/12\;\; \forall \epsilon\in (0, 1/2)\;.$$

\subsection{A restriction on partitions of unity}
\label{subsec-part}
\noindent

A subset $U \subset M$ is called {\it displaceable} if there
exists a Hamiltonian diffeomorphism $\psi$ of $M$ such that $\psi(U)
\cap \text{Closure}(U) = \emptyset$.

In the situation of Example \ref{ex-mfds} the quasi-state $\zeta$
has the following additional {\it vanishing property}: $\zeta(F) =
0$ provided $F$ has  a displaceable support\footnote{In
\cite{EP-qst}, the vanishing property, together with the
invariance under $\Symp_0(M,\omega)$, was included into the
definition of a symplectic quasi-state. Today we believe that they
should be considered as additional properties rather than axioms.}.

In this section we present a simple consequence of Theorem
\ref{thm-main-A} which, in particular, provides a restriction on
partitions of unity on $M$ subordinate to coverings of $M$ by
sufficiently small sets.

\begin{thm}\label{thm-part}
There exists a constant $K>0$, which depends only on the
symplectic manifold $(M,\omega)$, with the following property:
Given any $N$ functions $\rho_1,\ldots,\rho_N$ on $M$ with
displaceable supports so that $\sum_{i=1}^N \rho_i \geq 1,$ the
following inequality holds:
\begin{equation}
\label{eq-part} \max_{i,j} ||\{\rho_i,\rho_j\}|| \geq
\frac{K}{N^3}\;.
\end{equation}
\end{thm}

\begin{proof} Denote by $a: = \max_{i,j} ||\{\rho_i,\rho_j\}||$ the number in the left-hand side of \eqref{eq-part}.
For an integer $k \in [1;N]$ put $r_k = \rho_1+\ldots+\rho_k$.
Note that $\zeta(\rho_k) =0$ by the vanishing property. Thus
Theorem \ref{thm-main-A} yields
$$|\zeta(r_k)-\zeta(r_{k-1})| \leq
\sqrt{2C||\{r_{k-1},\rho_k\}||}\leq \sqrt{2Ca}\sqrt{k-1}\;.$$ Sum
up these inequalities for $k=2,\ldots,N$. Note that
$\zeta(r_N)\geq \zeta(1)=1$ in view of monotonicity axiom (ii) of
$\zeta$. Furthermore $\zeta(r_1)=\zeta(\rho_1)=0$. Hence we get
that
$$1 \leq \sqrt{2Ca}\sum_{k=2}^N \sqrt{k-1} \leq
\text{const}\cdot a^{\frac{1}{2}}N^{\frac{3}{2}}\;,$$ which yields
\eqref{eq-part}. \end{proof}

\medskip
\noindent {\sc A generalization:} Interestingly enough, a slightly
weaker version of Theorem \ref{thm-part} holds true for a much
more general class of closed symplectic manifolds $(M,\omega)$
than we considered before. For technical reasons we assume that
$M$ is {\bf rational}, i.e. the image of $\pi_2 (M)$ under the
cohomology class of $\omega$ is a discrete subgroup of $\R$.
Furthermore, we assume that $M$ is {\bf strongly semi-positive},
that is
\begin{equation} \label{eqn-strongly-semi-pos} 2-n\leq c_1 (A)<
0 \Longrightarrow \omega (A)\leq 0, \ {\rm for\ any}\ A\in \pi_2
(M),
\end{equation}
where $c_1$ stands for the 1st Chern class of $(M,\omega)$. For
instance, every symplectic $4$-manifold is strongly semi-positive.
We believe that eventually these assumptions will be omitted.

Fix a displaceable open subset $U\subset M$. A closed subset $X
\subset M$ is called {\it dominated by} $U$ if there exists a
Hamiltonian diffeomorphism $\psi$ of $M$ with $\psi(X) \subset U$.

\begin{thm}\label{thm-part-1}
There exists a constant $K>0$ which depends only on the symplectic
manifold $(M,\omega)$ and on the subset $U$ with the following
property: Given any $N$ functions $\rho_1,\ldots,\rho_N$ on $M$
whose supports are dominated by $U$ so that $\sum_{i=1}^N \rho_i
\geq 1\;,$ the following inequality holds:
\begin{equation}
\label{eq-part-prim} \max_{i,j} ||\{\rho_i,\rho_j\}|| \geq
\frac{K}{N^3}\;.
\end{equation}
\end{thm}

\medskip
\noindent The proof repeats the argument above with one
modification: a reference to Theorem \ref{thm-main-A} is replaced
by its weaker version, see Section \ref{subsec-proof-part-1} for
the details.

\subsection{Simultaneous measurability in
classical mechanics} \label{sect-measurement-general} \noindent
Symplectic quasi-states are classical analogues of quasi-states in
quantum mechanics. The latter appeared as an attempt to revise von
Neumann's notion of a quantum mechanical state as a {\it linear}
functional on the algebra of observables. A number of physicists
considered the equation $\xi(A+B) = \xi(A)+\xi(B)$, where $\xi$ is
a state and $A,B$ are observables, as lacking physical meaning
unless $A$ and $B$ commute: indeed, non-commuting observables are
not simultaneously measurable and hence the expression $\xi(A+B)$
is not well defined. This gave rise to the notion of quasi-state,
a non-linear functional which is linear on any subspace generated
by a pair of commuting observables (compare with axiom (iii) of a
symplectic quasi-state). We refer to \cite{EP-qst} for a detailed
historical account.

In view of this discussion, the existence of non-linear symplectic
quasi-states on classical observables (that is on functions on
symplectic manifolds) naturally leads us to the problem of
simultaneous measurability in classical mechanics. This problem
appears in physics literature (see e.g. books by Peres
\cite[Chapter 12-2]{Peres} and Holland \cite{Holland}) as a toy
example motivating the theory of quantum measurements. Below we
analyze simultaneous measurability in classical mechanics in the
framework of a measurement procedure, called the {\it pointer
model}. We shall show that $\Pi(F_1,F_2)$ gives a lower bound for
the error of simultaneous measurement of observables $F_1$ and
$F_2$.

Consider two observables $F_1,F_2 \in C^\infty (M)$. Let
\(\widehat M = M \times \R^{4}(p,q)\), $p = (p_1, p_2)$, $q= (q_1,
q_2)$, be the extended phase space equipped with the symplectic
form $\widehat \omega = \omega + dp\wedge dq$. The \(\R^4\) factor
corresponds to the measuring apparatus (the pointer), whereas
\(q\) is the quantity read from it. The coupling of the apparatus
to the system is carried out with the aid of the Hamiltonian
function \(p_1F_1(x)+p_2F_2(x)\). The Hamiltonian equations of
motion with the initial conditions $q(0)=0, p_1(0)=p_2(0) =
\epsilon$ and $x(0) = y$ are as follows:
\begin{align*}
\dot {q_i} &= F_i,\; i=1,2\\
\dot {p} &= 0\\
\dot x &= \epsilon\sgrad(F_1+F_2)\;.
\end{align*}
Denote by $g_t$ the Hamiltonian flow on $M$ generated by the
function $G=F_1+F_2$. Then $x(t) = g_{\epsilon t}y$. Let $T>0$ be
the duration of the measurement. By definition, the output of the
measurement procedure is a pair of functions $F'_i, \; i=1,2$, on
$M$ defined by the average displacement of the $q_i$-coordinate of
the pointer:
$$F'_i (y) = \frac{1}{T}(q_i(T) -q_i(0)) =\frac{1}{T}\int_0^T
F_i(x(t))dt =\frac{1}{T}\int_0^T F_i(g_{\epsilon t}y)dt \;.$$ Note
that for $\epsilon =0$ we have $F'_i =F_i$. This justifies the
above procedure as a measurement of $F_i$ and allows us to
interpret the number $\epsilon$ as {\it an imprecision of the
pointer}.

Define {\it the error of the measurement} as
\[\Delta(T,\epsilon,F_1,F_2) = \left\| F'_i - F_i \right\| .\]
Note that in our setting this quantity does not depend on $i \in
\{1;2\}$ since the sum $F_1+F_2$ is constant along the
trajectories of $g_t$.

\begin{thm}\label{thm-meas} For all $T,\epsilon > 0$ and $F_1,F_2
\in C^{\infty}(M)$
$$\Delta(T,\epsilon,F_1,F_2) \geq
\frac{1}{2}\Pi(F_1,F_2)-\sqrt{\frac{C}{T\epsilon}} \cdot
\sqrt{\min (\, \| F_1 -\langle F_1 \rangle \| ,  \| F_2 -\langle
F_2 \rangle \|\, ) }\;.$$
\end{thm}

\medskip
\noindent In particular, if $\Pi(F_1,F_2) \neq 0$ and the pointer
is not ideal, that is $\epsilon >0$, the error of the measurement
is bounded from below by $\Pi(F_1,F_2)/2$ when $T \to
+\infty$ uniformly in $\epsilon$. Theorem \ref{thm-meas} is proved
in Section \ref{subsec-proof-meas}.

\bigskip
\noindent {\sc ORGANIZATION OF THE PAPER:} Section \ref{sec-proofs}
contains proofs of Theorems \ref{thm-main-A}, \ref{thm-part-1} and
\ref{thm-meas}. Section \ref{subsec-vit} is a mock version of
Section \ref{subsec-Pois} where we revise lower bounds for the
Poisson bracket in terms of Hofer's geometry. Section \ref{sec-aux}
contains proofs of auxiliary facts on quasi-morphisms and
quasi-states used in the introduction. Finally, in Section
\ref{sec-prob} we discuss some open problems.

\section{Proofs}\label{sec-proofs}

\subsection{Proof of Theorem
\ref{thm-main-A}:}\label{subsec-proof-main}\noindent It suffices
to prove the theorem for functions $F$ and $G$ with zero mean. Let
\(f_t\) and \(g_t\) be the flows generated by \(F\) and \(G\). Put
\(H = F + G\) and \(K_t = F + G \circ f_t^{-1}\). Then \(K_t\) is
a normalized Hamiltonian generating the flow \(f_tg_t\) and so
$\phi_K = \phi_F\phi_G$. Observe that
\[\|H - K_t\| = \|G - G\circ f_t^{-1}\| = \|G\circ f_t - G\|\;.\]
Taking into account that
\[G(f_tx) - G(x) = \int_0^t\frac{d}{ds}G(f_sx)ds = -\int_0^t\{F,G\}(f_sx)ds \;,\]
we get that
\begin{equation}\label{eq-pr-main-1} \| H - K_t\| = \|G\circ f_t -
G\| \leq t\|\{F,G\}\|\;.
\end{equation}
From the stability property of \(\mu\) (see \eqref{eq-Lipschitz})
and inequality \eqref{eq-pr-main-1}, we get that
\begin{equation}\label{eq-vit-vsp}
|\mu(\phi_H) - \mu(\phi_K)| \leq \int_0^1 ||H-K_t||\;dt \leq
\|\{F,G\}\| \cdot \int_0^1 tdt = \frac{\|\{F,G\}\|}{2}\;.
\end{equation}
Combining this inequality with the fact that $C$ is the defect of
the quasi-morphism $\mu$ we obtain that
\begin{multline*}
\Pi(F,G) = |\mu(\phi_{F+G}) - \mu(\phi_F) - \mu(\phi_G)| \leq \\
\leq |\mu(\phi_H) - \mu(\phi_F\phi_G)| + |\mu(\phi_F) + \mu(\phi_G) - \mu(\phi_F\phi_G)| \leq \\
\leq |\mu(\phi_H) - \mu(\phi_K)| + C \leq \frac{\|\{F,G\}\|}{2} +
C.
\end{multline*}

Finally, let us balance this inequality. Let \(E > 0\) be any
number. Then
\[\Pi(EF,EG) \leq E^2\frac{\|\{F,G\}\|}{2} + C\;.\]
Since \(\Pi\) is homogeneous, after dividing both sides by \(E\)
we obtain
\[\Pi(F,G) \leq E\frac{\|\{F,G\}\|}{2} + \frac{C}{E}.\]
Choosing the optimal value $E = \sqrt{2C/||\{F,G\}||}$, we get
that
$$\Pi(F,G) \leq \sqrt{2C\|\{F,G\}\|}\,$$
as required. \Qed

\subsection{Proof of Theorem
\ref{thm-part-1}:}\label{subsec-proof-part-1}\noindent Let $\mu$
and $\zeta$ be the functionals introduced in
\cite[Section~7]{EP-qst}. These functionals satisfy a number of
properties which are weaker than the ones of a quasi-morphism and
of a quasi-state, but still these properties suffice to our
purposes. Put $\Pi(F,G)= |\zeta(F+G)-\zeta(F)-\zeta(G)|$. We claim
that there exists a constant $K_1>0$ so that
\begin{equation}\label{eq-part-1-vsp} \Pi(F,G) \leq K_1\sqrt{||\{F,G\}||}\end{equation}
for all $F,G \in C^{\infty}(M)$ so that the support of $G$ is
dominated by $U$. The proof repeats {\it verbatim} the proof of
Theorem \ref{thm-main-A} above. Now we repeat the proof of Theorem
\ref{thm-part} given in Section \ref{subsec-part} with one
modification: we replace the reference to  Theorem
\ref{thm-main-A} by the reference to \eqref{eq-part-1-vsp}. This
yields the desired result. \Qed

\subsection{Proof of Theorem \ref{thm-meas}:}\label{subsec-proof-meas}\noindent
It suffices to prove the theorem assuming that the functions $F_1$
and $F_2$ have zero mean. The proof is divided into several steps.

\medskip
\noindent{\sc Step 1:} One readily checks the following scaling
properties of the functional $\Delta(T,\epsilon,F_1,F_2)$:
\begin{equation}\label{eq-meas-1}
\Delta(T,\epsilon,F_1,F_2) = \Delta(\epsilon T,1,F_1,F_2)
\end{equation}
and
\begin{equation}\label{eq-meas-2}
\Delta(T,\epsilon,EF_1,EF_2) =
E\Delta(ET,\epsilon,F_1,F_2)\;\;\forall E>0\;.
\end{equation}
The proofs are straightforward and we omit them.

\medskip
\noindent{\sc Step 2:} Put
\[F' = \frac 1 T \int_0^T F_1\circ g_t \;dt \;.\]
For \(s \in [0,1]\)
\begin{align*}
\left\| F'\circ g_s - F' \right\| = \frac 1 T \left\| \int_0^T
F_1\circ g_{t+s}\; dt - \int_0^T F_1\circ g_t \;dt \right\|\\ =
\frac 1 T
\left\| \int_s^{T+s} F_1 \circ g_t\; dt - \int_0^T F_1 \circ g_t\; dt \right\|\\
= \frac 1 T \left\| \int_T^{T+s} F_1 \circ g_t\; dt - \int_0^s F_1
\circ g_t\; dt \right\|  \\ \leq \frac{2\|F_1\| s}{T}.
\end{align*}
Let \(f_s\) be the flow generated by \(F'\). Recall that
\(G=F_1+F_2\) generates the flow \(g_s\). Therefore \(-G + F'
\circ g_s\) generates the flow \(g_s^{-1}f_s\) whose time-one map
equals $\phi_G^{-1}\phi_{F'}$. Put \(K = F' - G\). Since $\mu$ is
stable (see formula \eqref{eq-Lipschitz}) and all our Hamiltonians
are normalized we obtain
\begin{align*}
|\mu(\phi_G^{-1} \phi_{F'})- \mu( \phi_K)| &\leq \int_0^1 \|-G + F' \circ g_s - (F' - G)\|ds\\
&= \int_0^1 \|F' \circ g_s - F' \|ds \leq \int_0^1
\frac{2\|F_1\|s}{T} ds \leq \frac{\|F_1\|}{T}.
\end{align*}

Since \(\mu\) is a homogeneous quasi-morphism with defect $C$, we
get
$$\Pi(F_1+F_2,-F') = \big|\mu(\phi_{F'}) - \mu(\phi_G) -
\mu(\phi_K)\big|$$ $$ \leq \big|\mu(\phi_G^{-1} \phi_{F'}) -
\mu(\phi_K)\big| + C \leq \frac{\| F_1\|} T + C.$$ Since $\Pi$ is
Lipschitz in both variables with respect to the uniform norm (see
formula \eqref{eq-Lip-Pi}),
$$|\Pi(F_1+F_2,-F_1) - \Pi(F_1+F_2,-F')| \leq
2||F_1-F'||=2\Delta(T,1,F_1,F_2)\;.$$ Observe that
$\Pi(F_1+F_2,-F_1) = \Pi(F_1,F_2)$. Combining this with the
previous inequality we get
\begin{equation}\label{eq-step2-2}
\Delta(T,1,F_1,F_2)\geq \frac{1}{2}\Pi(F_1,F_2) -\frac{\|
F_1\|}{2T} -\frac{C}{2}\;.
\end{equation}

\medskip
\noindent{\sc Step 3:} Using \eqref{eq-meas-1} and
\eqref{eq-step2-2} we get that
$$\Delta(T,\epsilon,F_1,F_2) \geq \frac{1}{2}\Pi(F_1,F_2) -\frac{||F_1||}{2T\epsilon}
-\frac{C}{2}\;.$$ We shall now balance this inequality. Let \(E >
0\) and $\tau = TE$. We have that
$$\Delta(T,\epsilon,EF_1,EF_2) \geq \frac{1}{2}\Pi(EF_1,EF_2) -\frac{E||F_1||}{2T\epsilon}
-\frac{C}{2}\;.$$ Recall that $\Pi(EF_1,EF_2) = E\Pi(F_1,F_2)$
and, in view of \eqref{eq-meas-2},$$ \Delta(T,\epsilon,EF_1,EF_2)
= E\Delta(ET,\epsilon,F_1,F_2)\;.$$ Substituting this into the
previous inequality, dividing by \(E\) and using $\tau = TE$ we obtain that
$$\Delta(\tau,\epsilon,F_1,F_2) \geq \frac{1}{2}\Pi(F_1,F_2)
-\frac{E||F_1||}{2\tau\epsilon} -\frac{C}{2E}\;.$$ This inequality
is true for every $\tau>0$ and $E>0$. Choosing the scaling factor
$E$ in the optimal way as $E = \sqrt{C\tau\epsilon/||F_1||}$ we
get that
$$\Delta(\tau,\epsilon,F_1,F_2) \geq
\frac{1}{2}\Pi(F_1,F_2)-\sqrt{\frac{||F_1||C}{\tau\epsilon}}\;.$$
Since $\Delta(T,\epsilon,F_1,F_2) = \left\| F'_1 - F_1 \right\| =
\left\| F'_2 - F_2 \right\|$ and $\Pi(F_1,F_2)$ are symmetric with
respect to $F_1, F_2$, switching $F_1$ and $F_2$ in the proof
above shows that
$$\Delta(\tau,\epsilon,F_1,F_2) \geq
\frac{1}{2}\Pi(F_1,F_2)-\sqrt{\frac{||F_2||C}{\tau\epsilon}}\;.$$
The last two inequalities immediately yield the needed result.\Qed

\section{Hofer's metric and the Poisson bracket}\label{subsec-vit}\noindent

In this section we revise the $C^0$-robustness of the Poisson
bracket (compare with Section \ref{subsec-Pois} and
\cite{Car-Vit}) from the viewpoint of Hofer's geometry on
$\Ham(M,\omega)$ (see e.g. \cite{Pbook} for an introduction). We
work on an arbitrary closed symplectic manifold $(M,\omega)$.
Denote $\osc F = \max_M F -\min_M F$. For a Hamiltonian $F$ on $M
\times [0;1]$ define $\psi_F \in \Ham(M,\omega)$ as the time one
map of the Hamiltonian flow generated by $F$. The group
$\Ham(M,\omega)$ carries a bi-invariant metric $\rho$ called the
Hofer metric which is defined as follows:
\[\rho(\psi_F,\psi_G) :=
\inf \int_0^1 \osc\, H_t \;dt,\] where \(H_t\) is a time-dependent
Hamiltonian generating \(\psi_F^{-1}\psi_G\) and the infimum is
taken over all such \(H_t\). One can easily see that
\begin{equation}\label{eq-hof-1}
\rho(\psi_F,\psi_G) \leq \int_0^1 \osc\, (F_t-G_t)\;dt\;
\end{equation}
for all  Hamiltonians $F$ and $G$.  For $a,b \in \Ham(M,\omega)$
write $[a,b]$ for the commutator $aba^{-1}b^{-1}$. Denote by $\id$
the unit element of $\Ham(M,\omega)$.

Take any $F,G \in C^{\infty}(M)$. Using inequality
\eqref{eq-hof-1} and arguing
 as in Section
\ref{subsec-proof-main} (compare with formula \eqref{eq-vit-vsp})
we get that $$\rho(\psi_F\psi_G,\psi_{F+G}) \leq \frac{\osc
\{F,G\}}{2}\;.$$ Switching $F$ and $G$ and using the bi-invariance
of $\rho$, we conclude that
\begin{equation}\label{eq-hof-2}
\rho(\id,[\psi_F,\psi_G]) \leq \osc \{F,G\}\;.
\end{equation}

Further,
$$[a,b][\alpha,\beta]^{-1} = (a\alpha^{-1})\cdot \alpha \Big{(}
(b\beta^{-1}) \cdot \beta \Big{(} (a^{-1}\alpha)\cdot
\alpha^{-1}(b^{-1}\beta)\alpha \Big{)}
\beta^{-1}\Big{)}\alpha^{-1}\;.$$ Together with the bi-invariance
of Hofer's metric and the triangle inequality this yields
\begin{equation}\label{eq-commut}
\rho([a,b],[\alpha,\beta])\leq 2 \rho(a,\alpha) + 2 \rho(b,\beta)\;.
\end{equation}

Take any $F',G' \in C^{\infty}(M)$.
It follows from \eqref{eq-commut} and
\eqref{eq-hof-1}
that
$$\rho([\psi_F,\psi_G],[\psi_{F'},\psi_{G'}]) \leq
 2 \osc\,(F - F')+ 2 \osc\,(G - G').$$
Applying inequality \eqref{eq-hof-2} to the pair of functions
$F^\prime$, $G^\prime$ and using the bi-invariance of $\rho$ we get
\begin{equation}\label{eq-hof-3}
\rho(\id,[\psi_F,\psi_G]) \leq \osc \{F', G'\} + 2 \osc\,(F - F') + 2 \osc\,(G - G').
\end{equation}
Taking into account that $\osc\, H \leq 2 ||H||$ we conclude that
\begin{equation}
\label{eq-hof-final}
\| \{ F', G'\} \| \geq  \frac{1}{2}\rho(\id,[\psi_F,\psi_G]) -
2 \| F - F^\prime\| - 2 \| G - G^\prime\|.
\end{equation}
Recalling
the notations of Section~\ref{subsec-Pois} we see that
\[
\epsilon_{max} (F, G) \geq  \frac{1}{4} \rho(\id,[\psi_F,\psi_G]), \ \
\Upsilon (F, G) \geq \frac{1}{2} \rho(\id,[\psi_F,\psi_G])
\]
 and moreover
\begin{equation} \label{eq-hof-4}
\Upsilon_{F,G} (\epsilon) \geq
\frac{1}{2} \rho(\id,[\psi_F,\psi_G])- 2 \epsilon\ \ \forall
\epsilon\in (0, \frac{1}{4} \rho(\id,[\psi_F,\psi_G])).
\end{equation}

Finally, assume that $F_m \to F,G_m \to G$ and $\{F_m,G_m\} \to 0$
in the uniform norm. It follows from \eqref{eq-hof-final} that
$\psi_F$ and $\psi_G$ commute. The same holds true for $\psi_{tF}$
and $\psi_{sG}$ with any $s,t \in \R$ and hence $\{F,G\}\equiv 0$. This
proves implication \eqref{eq-vit} which appears in \cite{Car-Vit}.

We refer to \cite{Muller-Oh,Oh,Vit} for other recent results
related to $C^0$-behavior of Hamiltonians.

\section{Auxiliary results on quasi-morphisms and
quasi-states}\label{sec-aux}

\subsection{Proof of monotonicity axiom (ii) of $\zeta$:}
\label{subseq-pr-eqlip}\noindent Given a stable homogeneous
quasi-morphism $\mu$ on $\cG$, let $\zeta$ be the functional
defined by \eqref{eq-qmm-qst}. We have to check that $\zeta(F)\leq
\zeta(G)$ for $F \leq G$. This follows immediately from the
inequality
\begin{equation}\label{eq-pr-eqlip-vsp}
\zeta(F)-\zeta(G) \leq \max(F-G)
\end{equation}
for all $F,G \in C^{\infty}(M)$. Indeed, put $F_0 =F-\langle F
\rangle$ and $G_0 = G-\langle G \rangle$ and observe that
$\zeta(F)=\zeta(F_0)+ \langle F \rangle $ and
$\zeta(G)=\zeta(G_0)+\langle G \rangle$. Since $F_0$ and $G_0$ are
normalized we can apply \eqref{eq-Lipschitz} and get that
$$\zeta(F_0)-\zeta(G_0)\leq \max(F_0-G_0)\;.$$
Taking into account that
$$\max (F_0-G_0) = \max(F-G) - \langle F-G \rangle\;,$$ we get
\eqref{eq-pr-eqlip-vsp}.  \Qed

\subsection{On the stability property
\eqref{eq-Lipschitz}}\label{subsec-stab}\noindent Here we show
that the homogeneous quasi-morphism $\mu$ constructed in
\cite{EP-qmm} is stable. First of all let us briefly recall the
construction of $\mu$. Let $QH$ be the even quantum homology
algebra of $M$. Recall that $\cG$ denotes the universal cover of
the group $\Ham (M,\omega)$.

Let
$$c: QH \times \cG \to \R$$
be {\it the spectral invariant} introduced by Y.-G. Oh
\cite{Oh-spectral,MS2}. Then
\begin{equation} \label{eq-mu}
\mu(f) = -\lim_{m \to \infty} \frac{c(e,f^m)}{m}\;
\end{equation}
for certain element $e \in QH$. (Recall that we normalize the symplectic
form so that the symplectic volume of $M$ is 1 and therefore $\mu$ is defined precisely as in
\cite{EP-qmm}).

It is known that
\begin{equation}\label{eq-Lipschitz-1}
\int_0^1 \min_M (F_t-G_t)\;dt \leq c(e,\phi_F)-c(e,\phi_G)
\leq \int_0^1 \max_M(F_t-G_t)\;dt
\end{equation}
for all normalized Hamiltonians  $F,G \in C^{\infty} (M \times
S^1)$, see \cite[formula (2.30)]{EP-qmm}.

Let $F,G$ be two normalized Hamiltonians on $M \times [0;1]$.
Without loss of generality assume that they are defined on $M
\times S^1$. This can be achieved by a suitable change of time in
the flows generated by $F$ and $G$ which alters the values of the
integrals in inequality \eqref{eq-Lipschitz-1} in an arbitrarily
small way.

Put $F_m(x,t) = mF(x,mt)$ and $G_m(x,t)= mG(x,mt)$ and note that
$\phi_F^m =\phi_{F_m}$ and $\phi_G^m =\phi_{G_m}$ for all $m \in
\N$. Applying \eqref{eq-Lipschitz-1} to $F_m$ and $G_m$ and
introducing the time variable $\tau = mt$ we get that
$$m\cdot \int_0^1\min_M (F_\tau-G_\tau)\;d\tau \leq c(e,\phi^m_F)-c(e,\phi^m_G)
\\\leq m\cdot \int_0^1 \max_M(F_\tau-G_\tau)\;d\tau\;.$$
Dividing by $m$ and passing to the limit as $m \to \infty$ we
conclude with the help of \eqref{eq-mu} that
$$ \int_0^1\min_M (F_{\tau}-G_{\tau})\;d\tau \leq \mu(\phi_G)-\mu(\phi_F)
\leq  \int_0^1 \max_M(F_{\tau}-G_{\tau})\;d\tau\;.$$
\Qed

\subsection{Estimating the defect for the
2-sphere}\label{sec-defect}\noindent Here we prove Proposition
\ref{prop-illustr}. We start with some preliminaries from
\cite{EP-qmm}: Consider the field  $k = \C [[s]$ of Laurent series
in one variable (the series are possibly infinite in the negative
direction but finite in the positive one). The quantum homology
algebra $QH$ of $\SP^2$ is the algebra of polynomials with
coefficients in $k$ in the variable $p$ modulo the ideal generated
by $p^2-s^{-1}$:
$$QH = k[p]/\{p^2 = s^{-1}\}\;.$$ The quasi-morphism $\mu$ is
defined by formula \eqref{eq-mu} above, where $e = 1 \in QH$.
Denote by $\id$ the unit element of $\cG$. We shall need the
following properties of the spectral invariant $c: QH \times \cG
\to \R$ which can be readily extracted from \cite{EP-qmm}:
\begin{itemize}
\item[{(i)}] $c(ps, \id) = 1$;
\item[{(ii)}] $c(ab,fg) \leq c(a,f)+c(b,g)$;
\item[{(iii)}] $c(1,g) = -c(p,g^{-1})$
\end{itemize}
for all $a,b \in QH$ and $ f,g \in \cG$.

\medskip
\noindent{\bf Proof of Proposition \ref{prop-illustr}:} We have to
show that the defect $C$ of $\mu$ satisfies $C \leq 2$. We claim
that
\begin{equation}\label{eq-defect}
c(1,f)+c(1,g)-1 \leq c(1,fg) \leq c(1,f)+c(1,g)
\end{equation}
for all $f,g \in \cG$. Indeed, the inequality on the right is an
immediate consequence of property (ii) above. To get the
inequality on the left we use (i), (ii) and (iii) and observe that
\begin{align*} c(1,f) = c(1\cdot p \cdot ps, fg \circ g^{-1} \circ
\id) \leq c(1,fg) +c(p,g^{-1})+c(ps,\id)\\= c(1,fg) -c(1,g) +
1\;.\end{align*} This proves \eqref{eq-defect}.

It follows from inequality \eqref{eq-defect} that
$$mc(1,h) -(m-1) \leq c(1,h^m) \leq mc(1,h)$$
and hence formula \eqref{eq-mu} yields
\begin{equation}\label{eq-defect-1}
c(1,h)-1 \leq -\mu(h) \leq c(1,h)
\end{equation}
for all $h \in \cG$.

Take any $f,g \in \cG$. Substituting consecutively $h=f$, $h=g$
and $h=fg$ into inequality \eqref{eq-defect-1} and using
\eqref{eq-defect} we readily get that
$$-2\leq \mu(fg)-\mu(f)-\mu(g) \leq 2\;,$$
which completes the proof. \Qed

\section{Discussion and open problems}\label{sec-prob}

In Section \ref{subsec-Pois} we have defined the functional
$$\Upsilon(F,G) = \liminf_{F',G'\stackrel{C^0}{\longrightarrow}F,G}
||\{F',G'\}||\;.$$ The main open problem concerning this functional
is as follows:

\begin{question}\label{quest-1}
Is it true that we always have \(\Upsilon(F,G) = \|\{F,G\}\|\)?
\end{question}

\medskip
\noindent In a recent work by one of the authors \cite{Zap} this
question is answered in the positive for two-dimensional symplectic
manifolds using methods of the topology of surfaces.

In the case when $M$ is the sphere $\SP^2$, the symplectic
quasi-state $\zeta$, and hence the functional $\Pi$, are defined in
elementary terms (see Example \ref{ex-qst-sphere} above). Thus the
{\it formulations} of Theorems \ref{thm-main-A} and \ref{thm-meas}
are ``soft''. However our proofs of these theorems are
``hard''\footnote{The terms ``soft'' and ``hard'' are understood
here in the sense of Gromov \cite{Gromov-ICM}.}: they use in a
crucial way the fact that $\zeta$ is induced by a stable homogeneous
quasi-morphism on $\cG$ which is defined by means of Floer homology.
This somewhat paradoxical situation was partially resolved in
\cite{Zap}: For the sphere, Theorem \ref{thm-main-A} is proved in
\cite{Zap} by ``soft'' methods. Moreover, an analogous theorem is
proved with the sphere replaced by an arbitrary closed symplectic
surface \(\Sigma\), and with the quasi-state of Theorem
\ref{thm-main-A} replaced by an arbitrary simple quasi-state (that
is a quasi-state which is \emph{multiplicative} on each singly
generated closed subalgebra of \(C(\Sigma)\)). It is still unclear
to us whether Theorem \ref{thm-meas} on the error of a simultaneous
measurement of non-commuting observables admits a ``soft" proof in
the case of $\SP^2$.

At the moment the authors believe that in higher dimensions, the
$C^0$-robustness of the Poisson bracket is a ``hard'' phenomenon and
thus the Floer-theoretical methods used above are adequate in this
context.

Another interesting circle of problems is  related to the sharpness of
various estimates obtained in the present work. For instance,
Zapolsky's result mentioned in the beginning of this section shows
that inequality  \eqref{eq-illustr-1} stating that
$\Upsilon_{x^2,y^2}(\epsilon) \geq (1-2\epsilon)^2/4$ is not
asymptotically (as $\epsilon \to 0$) sharp: indeed one readily
computes that $||\{x^2,y^2\}|| \approx 9.7 > 0.25$. It would be
interesting to understand whether the approach of \cite{Zap} gives
rise to sharp lower bounds for $\Upsilon_{x^2,y^2}(\epsilon)$ on
surfaces.

\medskip
\noindent \begin{question} What is the sharp value of the constant
$C$ in Theorem \ref{thm-main-A} stating that $ \Pi(F,G) \leq
\sqrt{2C \cdot ||\{F,G\}||}$?
\end{question}

\medskip
\noindent The answer is so far unknown even for the case of the
quasi-state $\zeta$ from Example \ref{ex-qst-sphere} on the
2-sphere. In this case Floer-theoretical Proposition
\ref{prop-illustr} above yields $C \leq 2$, while the topological
argument from \cite{Zap} improves this to $C \leq 1/2$.

\medskip
\noindent
\begin{question}\label{quest-3}
Can one improve (asymptotically in $N$ as $N \to \infty$)
the bound
$$\max_{i,j} ||\{\rho_i,\rho_j\}|| \geq \text{const}/N^3$$ for the
partitions of unity given by Theorem \ref{thm-part}?
\end{question}

\medskip
\noindent  We claim that the asymptotic behavior of the right hand
side cannot be made better than $\sim N^{-2}$. Indeed, fix any
partition of unity $\rho'_1,...,\rho'_d$ subordinate to a
covering by displaceable subsets. For any $m \in \N$ consider
$N=md$ functions $\rho_i$ where $\rho_i = m^{-1}\cdot \rho'_j$ for
$i = j\; ({\rm mod}\; d).$ In other words, we take each function
$\rho'_j/m$ with multiplicity $m$. Of course, we again get a
partition of unity subordinate to the covering by displaceable
subsets. When $m \to \infty$, the left hand side of the inequality
above is of the order $\sim m^{-2}= d^2 N^{-2}$, and the claim
follows. At the moment, we do not know an answer to Question
\ref{quest-3} even in dimension two.

\bibliographystyle{alpha}

\bigskip

\noindent
\begin{tabular}{l}
Michael Entov \\
Department of Mathematics\\
Technion - Israel Institute of Technology \\
Haifa 32000, Israel\\
entov@math.technion.ac.il \\
\end{tabular}

\bigskip
\noindent
\begin{tabular}{l}
Leonid Polterovich \\
School of Mathematical Sciences \\
Tel Aviv University \\
Tel Aviv 69978, Israel \\
polterov@post.tau.ac.il \\
\end{tabular}

\bigskip
\noindent
\begin{tabular}{l}
Frol Zapolsky \\
School of Mathematical Sciences \\
Tel Aviv University \\
Tel Aviv 69978, Israel \\
zapolsky@post.tau.ac.il \\
\end{tabular}

\end{document}